\documentclass[a4paper,10pt]{article}
\usepackage{amsmath,amssymb,amsfonts}
\usepackage{epsfig,cite}
\usepackage{indentfirst}
\usepackage{algorithm2e}

\makeatletter
\def\@makefnmark{\hbox{\@textsuperscript{\bf\normalfont(\@thefnmark}}}
\makeatother

\usepackage [cp1251]{inputenc}
\usepackage[english]{babel}

\newtheorem{theorem}{Theorem}[section]
\numberwithin{equation}{section}
\numberwithin{theorem}{section}

\title{The FD-method for solving Sturm-Liouville problems with special singular differential operator}
\author{Volodymyr L. Makarov, Denis V. Dragunov, Yaroslav V. Klimenko }
\begin{document}
\maketitle
\large
\selectlanguage{english}
\begin{abstract}In the paper we describe a superexponentially convergent numerical-analytical method for solving the eigenvalue problem for the some class of singular differential operators with boundary conditions. The method (FD-method) was firstly proposed by V. L. Makarov and successfully combines the benefits of using the {\it coefficient approximation methods} (CAM) and the homotopy approach.  The sufficient conditions which provides convergence of the proposed method are stated and rigorously substantiated. The algorithm for the software implementation of the proposed method is described too. A lot of numerical examples are included in the paper. The examples confirm the theoretical conclusions. We also have made the comparison between the results obtained by FD-method and results obtained by  the powerful software package for solving Sturm-Liouville problems --- SLEIGN2.
\end{abstract}


\section{Introduction}
There is a great scope of the numerical methods for solving eigenvalue problems for differential operators of the second order (see \cite{pryce}). The known numerical techniques for eigenvalue problems can be divided into two groups: methods based on the direct approximation of the solution of differential equation and methods based on approximation of the coefficients of differential equation (see \cite{Makarov_Rossokhata_Review}).

The methods of the first group such as finite difference, finite elements and spectral are extensively treated through both theoretical investigations and well developed soft. The main idea of this approach is replacement of the eigenfunctions by piecewise polynomial functions, which results in the approximation of a differential equation by a system of linear algebraic equations. Because of the different nature of the original differential operator and approximating algebraic operator the numerical solution comes close to the exact one only for the low-indexed eigenvalues and corresponding eigenfunctions, that is defined by the size of the mesh step. Hence, such numerical techniques are effective to find the low-indexed eigenvalues, but not effective to find the high-indexed eigenvalues \cite{Makarov_Rossokhata_Review}.

The methods of the second group are based on the general idea of the {\it coefficient approximation methods} (hereinafter referred to as CAM), also known as the Pruess methods, --- to replace the coefficient functions of the problem by simpler approximation and then to solve the approximating problem (see \cite{pryce}). Firstly this idea was proposed by Kryloff and Bogolioubov (1926) in \cite{kryloff_bogoliubov} for the case of piecewise constant approximation.  S. Pruess provided a thorough convergence and error analysis of such methods when piecewise polynomial approximation is used. Except for Gordon, who uses linear functions, the earlier  references  all confine themselves to piecewise constant approximations, and this is the most practical choice: otherwise the approximating problem may be no easier to solve numerically than is the original. It is well known (see \cite{pryce}) that if the piecewise constant approximation with the maximum step size $h$ is used, the error of the CAM is estimated by $O(h).$ Hence, when we using the CAM, a crucial problem arise: to halve the error we have to double the computational time. This problem in many cases can be solved by using {\it the functional-discrete} method (hereinafter referred to as the {\it FD-method}) which is the essential generalization of the CAM for solving Sturm-Liouville problems of different types. The algorithm of FD-method was firstly described and justified for linear Sturm-Liouville problem in \cite{makarov1}. Generally speaking, it consists of two main steps: the first one is to find the initial approximation by applying the CAM and the second one is to calculate a sufficient number of corrections to achieve needed accuracy. The second step is substantiated by the {\it homotopy method} (also known as { \it perturbation method}, {\it continuation method} or {\it successive loading method}, see, for example, \cite{Homotopy_NATO}). From that point of view, FD-method is the synthesis of the CAM and homotopy method.

In the present paper we consider the following singular Sturm-Liouville problem
\begin{equation}\label{intro_exp_1}
    \frac{d}{d x}\left[\left(1-x^{2}\right)\frac{d u(x)}{d x}\right]-q(x)u(x)=\lambda u(x),\quad x\in \left(-1,1\right),
\end{equation}
\begin{equation}\label{intro_exp_2}
    \lim\limits_{x\rightarrow \pm 1}\left(1-x^{2}\right)\frac{d u(x)}{d x}=0.
\end{equation}

% Applying the FD-method, we looking for the eigenvalues and %eigenfunction in the form of the series %$\lambda=\sum\limits_{i=0}^{\infty}\lambda^{(i)},$ %$u(x)=\sum\limits_{i=0}^{\infty}u^{(i)}(x)$ respectively (see %\cite{makarov1}).

Since the problem \eqref{intro_exp_1},\eqref{intro_exp_2} has singularities  at the points $\pm 1,$ we can't directly apply the algorithm of the FD-method, proposed in \cite{makarov1}. However, in \cite{Klimenko_Ya}  the modified FD-method technique for the problems with such singularities have been presented. Let us state briefly the algorithm of the FD-method described in \cite{Klimenko_Ya}.
According to the FD-method we approximate the eigenvalues and eigenfunctions of the problem \eqref{intro_exp_1}, \eqref{intro_exp_2} by the following truncated series
\begin{equation}\label{intro_exp_3}
    \stackrel{m}{\lambda}_{n}=\sum\limits_{j=0}^{m}\lambda_{n}^{(j)}, \quad \stackrel{m}{u}_{n}(x)=\sum\limits_{j=0}^{m}u_{n}^{(j)}(x),
\end{equation}
where the unknown summands $\lambda_{n}^{j}$ and  $u_{n}^{(j)}(x),$ for $i=1,2,\ldots, m\footnote{The positive integer number $m$ is called the {\it rank} of FD-method.}$ are the solutions to the following recurrence problems:
\begin{equation}\label{intro_exp_4}
    \frac{d}{dx}\left[\left(1-x^{2}\right)\frac{d u_{n}^{(j)}(x)}{d x}\right]+\lambda_{n}^{(i)}u_{n}^{(j)}(x)=-\sum\limits_{s=0}^{j-1}\lambda_{n}^{(j-s)}u_{n}^{(s)}(x)+q(x)u_{n}^{(j-1)}(x)=F^{(j)}_{n}(x),
\end{equation}
$x\in (-1,1)$ with the boundary conditions
\begin{equation}\label{intro_exp_5}
    \lim\limits_{x\rightarrow \pm 1}\left[\left(1-x^{2}\right)\frac{d u_{n}^{(j)}(x)}{d x}\right]=0,
\end{equation}
and the orthogonality conditions
\begin{equation}\label{intro_exp_6}
    \int\limits_{-1}^{1}u^{(j)}_{n}(x)u_{n}^{(0)}(x)d x=0, \quad j=1,2,\ldots, m.
\end{equation}

The initial data $\lambda_{n}^{(0)},$ $u_{n}^{(0)}(x)$ which is needed to begin the recurrence process can be found as the solution to the following eigenvalue problem, called {\it the basic problem}:
\begin{equation}\label{intro_exp_7}
    \frac{d}{d x}\left[\left(1-x^{2}\right)\frac{d u^{(0)}_{n}(x)}{d x}\right]+\lambda_{n}^{(0)}u^{(0)}_{n}(x)=0,\quad x\in\left(-1,1\right),
\end{equation}
\begin{equation}\label{intro_exp_8}
    \lim\limits_{x\rightarrow \pm 1}\left[\left(1-x^{2}\right)\frac{d u^{(0)}_{n}(x)}{d x}\right]=0.
\end{equation}
The solutions to problem \eqref{intro_exp_7}, \eqref{intro_exp_8} can be expressed in the following form
\begin{equation}\label{intro_exp_9}
    \lambda_{n}^{(0)}=n\left(n+1\right),\quad u_{n}^{(0)}(x)=\sqrt{\frac{2n+1}{2}}P_{n}(x),
\end{equation}
where $P_{n}(x)$ denotes the Legendre polynomial of the first kind of order $n\in \mathbb{N}$ (see \cite{Bateman_1}).
The following statement about the convergence of the FD-method, described above, was proved in \cite{Klimenko_Ya}.

\begin{theorem}\label{teor_1}
Let the function $q(x)$ belongs to the $Q^{0}\left[-1,1 \right],$ --- the space of piecewise continuous functions defined on $\left[-1, 1\right].$ If there exists a nonnegative integer $n_{0},$ such that
\begin{equation}\label{intro_exp_10}
    q_{n_{0}}=\frac{2\left\|q(x)\right\|_{\infty, \left[-1,1\right]}}{n_{0}}<1,\quad n_{0}\geq 1,\quad q_{0}=2\left\|q(x)\right\|_{\infty, \left[-1,1\right]}<1,
\end{equation}
then the FD-method \eqref{intro_exp_4} --\eqref{intro_exp_9} has a superexponential convergence rate for all numbers $n\geq n_{0}:$
\begin{equation}\label{intro_exp_11}
    \left\|u_{n}(x)-\stackrel{m}{u}_{n}(x)\right\|\leq \frac{\left(q_{n}\right)^{m+1}}{1-q_{n}}\alpha_{m+1},\quad \left|\lambda_{n}-\stackrel{m}{\lambda}_{n}\right|\leq \frac{\left(q_{n}\right)^{m}}{1-q_{n}}\left\|q(x)\right\|_{\infty}\alpha_{m},
\end{equation}
$n=n_{0}, n_{0}+1,\ldots,$ where
$$\alpha_{m+1}=2\frac{(2m-1)!!}{(2m+2)!!}\leq \frac{1}{(m+1)\sqrt{\pi m}},\quad \left\|q(x)\right\|_{\infty, \left[-1,1\right]}=\max\limits_{x\in\left[-1, 1\right]}\left|q(x)\right|.$$
\end{theorem}
Theorem \ref{teor_1} was proved in \cite{Klimenko_Ya} for the case when $q(x)=q(-x)$ but, using the same technique, it can be proved for the case of an arbitrary function from $Q^{0}\left[-1, 1\right]$ without considerable difficulties.

If we try to apply the FD-method \eqref{intro_exp_4} --\eqref{intro_exp_9} to approximate eigensolution $\lambda_{n}, u_{n}(x)$ with $n=0,1,\ldots, n_{0}-1$ we might find it divergent. The reason of that is the lack of approximation for $q(x)$ in equation \eqref{intro_exp_7} of the basic problem. To overcome such difficulties we have to use the general scheme of the FD-method which uses the approximation of the function $q(x)$ on $\left[-1, 1\right]$ by the piecewise constant function $\bar{q}(x).$

\section{The general algorithm of the FD-method for solving the Sturm-Liouville problem with Legendre differential operator: theoretical justification}

To make the general algorithm of the FD-method more easy for understanding, let us consider the following general problem
\begin{equation}\label{gen_alg_1}
    \frac{\partial}{\partial x}\left[\left(1-x^{2}\right)\frac{\partial u(x, t)}{\partial x}\right]+\left[\lambda(t)-\bar{q}(x)-t\left(q(x)-\bar{q}(x)\right)\right]u(x,t)=0, \quad x\in\left(-1,1\right),\; t\in\left[0,1\right],
\end{equation}
\begin{equation}\label{gen_alg_2}
    \lim\limits_{x\rightarrow \pm 1}\left(1-x^{2}\right)\frac{\partial u(x,t)}{\partial x}=0,\quad \forall t\in\left[0,1\right].
\end{equation}
To define the function $\bar{q}(x),$ we have to fix some mesh on the interval $[-1,1]:$
\begin{equation}\label{gen_alg_mesh}
    \bar{\omega}=\left\{-1=x_{0}<x_{1}<\ldots < x_{N-1}<x_{N}=1\right\},\quad h=\max\limits_{i=1,2,\ldots, N}\left(x_{i}-x_{i-1}\right).
\end{equation}
We also require that all points of discontinuity of the function $q(x)$ on the interval $[-1, 1]$ are contained in the set $\bar{\omega}.$

There are many possibilities to define  $\bar{q}(x),$ for example, $$\bar{q}(x)=\bar{q}(x, \bar{\omega})=q\left(\frac{1}{2}\left(x_{i-1}+x_{i}\right)\right),\quad x \in \left[x_{i-1}, x_{i}\right),\quad i=1,2,\ldots, N,$$
or
\begin{equation}\label{gen_alg_q}
   \bar{q}(x)=\bar{q}(x, \bar{\omega})=\frac{1}{2}\left(q(x_{i-1})+q(x_{i})\right),\quad x \in \left[x_{i-1}, x_{i}\right), \quad i=1,2,\ldots, N,
\end{equation}
but, in any case, the piecewise constant  function $\bar{q}(x,\bar{\omega})$ must have the following property
$$\lim\limits_{h\rightarrow 0}\left(\max\limits_{x\in\left[-1, 1\right]}\left|\bar{q}(x,\bar{\omega})-q(x)\right|\right)=0.$$
It is easy to see that if we set $t=1,$ then problem \eqref{gen_alg_1}, \eqref{gen_alg_2} would transform to the problem \eqref{intro_exp_1}, \eqref{intro_exp_2}.
Assume that the solution to the problem \eqref{gen_alg_1}, \eqref{gen_alg_2} can be expressed in the series form
\begin{equation}\label{gen_alg_3}
     \lambda_{n}=\sum\limits_{j=0}^{\infty}t^{(j)}\lambda_{n}^{(j)}, \quad u_{n}(x,t)=\sum\limits_{j=0}^{\infty}t^{j}u_{n}^{(j)}(x),
\end{equation}
and
\begin{equation}\label{gen_alg_4}
   \frac{\partial u_{n}(x,t)}{\partial x}=\sum\limits_{j=0}^{\infty}t^{j}\frac{\partial u_{n}^{(j)}(x)}{\partial x},\quad \frac{\partial^{2} u_{n}(x,t)}{\partial x^{2}}=\sum\limits_{j=0}^{\infty}t^{j}\frac{\partial^{2} u_{n}^{(j)}(x)}{\partial x^{2}},\quad \forall x\in\left(-1,1\right),\; \forall t\in\left[0,1\right].
\end{equation}
The unknown summands of the series \eqref{gen_alg_3} can be found as the solutions to the recurrence problems ($j=0,1,\ldots, m$):
\begin{equation}\label{gen_alg_5}
    \frac{d}{d x}\left[\left(1-x^{2}\right)\frac{d u_{n}^{(j)}(x)}{d x}\right]+\left[\lambda_{n}^{(0)}-\bar{q}(x)\right]u^{(j)}_{n}(x)=
\end{equation}
$$=-\sum\limits_{s=0}^{j-1}\lambda_{n}^{(j-s)}u_{n}^{(s)}(x)+\left[q(x)-\bar{q}(x)\right]u_{n}^{(j-1)}(x)\equiv F^{(j)}_{n}(x), \quad x\in (-1,1),$$
with the boundary conditions
\begin{equation}\label{gen_alg_6}
    \lim\limits_{x\rightarrow \pm 1}\left[\left(1-x^{2}\right)\frac{d u^{(j)}_{n}(x)}{d x}\right]=0,
\end{equation}
matching conditions
\begin{equation}\label{gen_alg_7_mk}
    \left[u_{n}^{(j)}(x)\right]_{x=x_{i}}=0,\quad \left[\frac{d u_{n}^{(j)}(x)}{d x}\right]_{x=x_{i}}=0, i=1,2,\ldots, N-1.
\end{equation}
and additional requirement
\begin{equation}\label{gen_alg_7}
     \int\limits_{-1}^{1}u^{(j)}_{n}(x)u^{(0)}_{n}(x)d x=\delta_{0,j}=\left\{
    \begin{array}{cc}
    1, & \mbox{if}\; j=0, \\
    0, & \mbox{if}\; j\neq 0.\\
    \end{array}
    \right.
\end{equation}
Here square brackets denotes the jump of the function at the point $x=x_{i}.$

Among the problems \eqref{gen_alg_5} -- \eqref{gen_alg_7_mk} it is worth to emphasize the first one ($j=0$):
\begin{equation}\label{gen_alg_8}
    \frac{d}{d x}\left[(1-x^{2})\frac{d u^{(0)}_{n}(x)}{d x}\right]+\left[\lambda_{n}^{(0)}-\bar{q}(x)\right]u^{(0)}_{n}(x)=0,\quad x\in (-1,1),
\end{equation}
\begin{equation}\label{gen_alg_9}
    \lim\limits_{x\rightarrow \pm 1}\left[\left(1-x^{2}\right)\frac{d u^{(0)}_{n}(x)}{d x}\right]=0,\quad \int\limits_{-1}^{1}\left(u_{n}^{(0)}(x)\right)^{2} d x=1,
\end{equation}
\begin{equation}\label{gen_alg_9_mk_bp}
    \left[u_{n}^{(0)}(x)\right]_{x=x_{i}}=0,\quad \left[\frac{d u_{n}^{(0)}(x)}{d x}\right]_{x=x_{i}}=0.
\end{equation}
As it was mentioned above, the problem \eqref{gen_alg_8} -- \eqref{gen_alg_9_mk_bp} is called {\it the basic problem}. It is easy to make sure that the differential operator
\begin{equation}\label{gen_alg_10}
    L(\cdot)=\frac{d}{d x}\left[\left(1-x^{2}\right)\frac{d (\cdot)}{d x}\right]-\bar{q}(x)(\cdot)
\end{equation}
is self-adjoint in the Sobolev space $$\overline{W}^{2,1}\left(-1,1\right)=\left\{f(x)\in W^{2,1}\left(-1,1\right)\colon \lim\limits_{x\rightarrow \pm 1}\left(1-x^{2}\right)f^{\prime}(x)=0\right\}.$$ This fact implies that the eigenfunctions $u_{n}^{(0)}(x),$ $n=0,1,\ldots$ of the problem \eqref{gen_alg_8} -- \eqref{gen_alg_9_mk_bp} form the complete orthogonal system in $\overline{W}^{2,1}\left(-1,1\right)$  and the corresponding  eigenvalues $\lambda_{n}^{(0)},$ ($\lambda_{i}^{(0)}<\lambda_{j}^{(0)},\;\mbox{when}\; i<j,\; i,j=0,1,2\ldots$) are simple.

%The following estimate holds true (see \cite{Ukhanev_Mak})
%\begin{equation}\label{gen_alg_11}
%    \left|\lambda_{n}^{(0)}-\lambda_{n}\right|\leq %\left\|q(x)-\bar{q}(x)\right\|_{\infty, \left[-1,1\right]}.
%\end{equation}
Suppose that the basic problem is solved and a pare $\lambda_{n}^{(0)}, u_{n}^{(0)}(x)$ denotes its arbitrary eigensolution.

It is well known that the problems \eqref{gen_alg_5} -- \eqref{gen_alg_7_mk} for $j=1,2,\ldots, m$ can be solved if and only if the following equality holds
\begin{equation}\label{gen_alg_12}
    \int\limits_{-1}^{1}u^{(0)}_{n}(x)F^{(j)}_{n}(x) dx=0.
\end{equation}
Equality \eqref{gen_alg_12} together with \eqref{gen_alg_5} yield us the formula for computing $\lambda_{n}^{(j)}:$
\begin{equation}\label{gen_alg_13}
    \lambda_{n}^{(j)}=-\sum\limits_{s=1}^{j-1}\lambda_{n}^{(j-s)}\int\limits_{-1}^{1}u^{(s)}_{n}(x)u_{n}^{(0)}(x) d x +\int\limits_{-1}^{1}\left(q(x)-\bar{q}(x)\right)u_{n}^{(j-1)}(x)u_{n}^{(0)}(x)d x.
\end{equation}
Taking into account \eqref{gen_alg_7} we can simplify formula \eqref{gen_alg_13} in the following way
\begin{equation}\label{gen_alg_14}
    \lambda_{n}^{(j)}=\int\limits_{-1}^{1}\left(q(x)-\bar{q}(x)\right)u_{n}^{(j-1)}(x)u_{n}^{(0)}(x)d x.
\end{equation}
The function $u_{n}^{(j)}(x)$ can be represented as a Fourier series
\begin{equation}\label{gen_alg_15}
   u_{n}^{(j)}(x)=\sum\limits_{s=0}^{\infty}u^{(0)}_{s}(x)\int\limits_{-1}^{1}u^{(j)}_{n}(x)u_{s}^{(0)}(x) d x.
\end{equation}
To find the Fourier coefficients let us multiply both sides of the equation \eqref{gen_alg_5} by $u_{s}^{(0)}(x)$ and then integrate them on the interval $[-1,1]:$
\begin{equation}\label{gen_alg_15}
    \int\limits_{-1}^{1}u_{s}^{(0)}(x)\frac{d}{d x}\left[\left(1-x^{2}\right)\frac{d u_{n}^{(j)}(x)}{d x}\right]d x+\int\limits_{-1}^{1}u_{s}^{(0)}(x)\left[\lambda_{n}^{(0)}-\bar{q}(x)\right]u_{n}^{(j)}(x)d x=
\end{equation}
$$=\int\limits_{-1}^{1}u_{s}^{(0)}(x)F^{(j)}_{n}(x) d x.$$
Taking into account boundary conditions \eqref{gen_alg_6} and using the integration by parts (two times), from \eqref{gen_alg_15} we obtain
\begin{equation}\label{gen_alg_16}
    \int\limits_{-1}^{1}\left(u_{n}^{(j)}(x)\frac{d}{d x}\left[\left(1-x^{2}\right)\frac{d u_{s}^{(0)}(x)}{d x}\right]-\bar{q}(x)u_{s}^{(j)}(x)u_{s}^{(0)}(x)\right)d x+\lambda_{n}^{(0)}\int\limits_{-1}^{1}u_{n}^{(j)}(x)u_{s}^{(0)}(x)dx=
\end{equation}
$$=\int\limits_{-1}^{1}u_{s}^{(0)}(x)F^{(j)}_{n}(x) d x.$$
 And \eqref{gen_alg_16} immediately lead us to the equality
 \begin{equation}\label{gen_alg_17}
    \int\limits_{-1}^{1}u_{n}^{(j)}(x)u_{s}^{(0)}(x)d x=\frac{\int\limits_{-1}^{1}F^{(j)}_{n}(x)u_{s}^{(0)}(x) d x}{\lambda_{n}^{(0)}-\lambda_{s}^{(0)}}
 \end{equation}
 Using \eqref{gen_alg_17} we can express the formula \eqref{gen_alg_15} in the following form
 \begin{equation}\label{gen_alg_18}
    u_{n}^{(j)}(x)=\sum\limits_{\substack{p=0\\p\neq n}}^{\infty}\frac{\int\limits_{-1}^{1}F^{(j)}_{n}(x)u_{p}^{(0)}(x) d x}{\lambda_{n}^{(0)}-\lambda_{p}^{(0)}}u_{p}^{(0)}(x).
 \end{equation}
 Formulas \eqref{gen_alg_14} and \eqref{gen_alg_18} yields the estimations
 \begin{equation}\label{gen_alg_19}
    \left|\lambda_{n}^{(j)}\right|\leq \left\|q(x)-\bar{q}(x)\right\|_{\infty, \left[-1, 1\right]}\left\|u_{n}^{(j-1)}(x)\right\|,
 \end{equation}
 and
 \begin{equation}\label{gen_alg_20}
    \left\|u_{n}^{(j)}(x)\right\|\leq M_{n}\left\|F^{(j)}_{n}(x)\right\|=M_{n}\left\{\left\|\sum\limits_{s=1}^{j-1}\lambda_{n}^{j-s}u_{n}^{(s)}(x)-\left(q(x)-\bar{q}(x)\right)u_{n}^{(j-1)}(x)\right\|^{2}-\left(\lambda^{j+1}\right)^{2}\right\}^{\frac{1}{2}}\leq
 \end{equation}
 $$\leq M_{n}\left\{\sum\limits_{s=1}^{j-1}\left|\lambda_{n}^{(j-s)}\right|\left\|u_{n}^{(s)}(x)\right\|+\left\|q(x)-\bar{q}(x)\right\|_{\infty, [-1,1]}\left\|u_{n}^{(j-1)}(x)\right\|\right\}\leq$$
 $$\leq M_{n}\left\|q(x)-\bar{q}(x)\right\|_{\infty, [-1,1]}\left\{\sum\limits_{s=0}^{j-1}\left\|u_{n}^{(j-1-s)}\right\|\left\|u_{n}^{(s)}(x)\right\|\right\},$$
 respectively, where
 \begin{equation}\label{gen_alg_20_Mn}
    M_{n}=\max\left\{\left(\lambda_{n}^{(0)}-\lambda_{n-1}^{(0)}\right)^{-1}, \left(\lambda_{n+1}^{(0)}-\lambda_{n}^{(0)}\right)^{-1}\right\}, \quad n=1,2,\ldots,
 \end{equation}
 $$M_{0}=\left(\lambda_{1}^{(0)}-\lambda_{0}^{(0)}\right)^{-1},$$ and $\left\|\cdot\right\|$ denotes the common $L_{2}$-norm, $\left\|f(x)\right\|=\sqrt{\int\limits_{-1}^{1} f^{2}(x)d x}.$

 Multiplying both sides of the inequality \eqref{gen_alg_20} by $\left(M_{n}\left\|q(x)-\bar{q}(x)\right\|_{\infty, [-1,1]}\right)^{-j}$ and using estimate \eqref{gen_alg_19}, we obtain
 \begin{equation}\label{gen_alg_21}
    \frac{\left\|u_{n}^{(j)}(x)\right\|}{\left(M_{n}\left\|q(x)-\bar{q}(x)\right\|_{\infty, [-1,1]}\right)^{j}}\leq\sum\limits_{s=0}^{j-1}\frac{\left\|u_{n}^{(j-s-1)}(x)\right\|\left\|u_{n}^{(s)}(x)\right\|}{\left(M_{n}\left\|q(x)-\bar{q}(x)\right\|_{\infty, [-1,1]}\right)^{j-1}}.
 \end{equation}
Using the notation
\begin{equation}\label{gen_alg_22}
    v_{j}=\left\|u_{n}^{(j)}(x)\right\|\left(M_{n}\left\|q(x)-\bar{q}(x)\right\|_{\infty, [-1,1]}\right)^{-j}=\left\|u_{n}^{(j)}(x)\right\|\bar{r}_{n}^{-j},
\end{equation}
we can represent inequality \eqref{gen_alg_21} in the following form
\begin{equation}\label{gen_alg_23}
    v_{j}\leq \sum\limits_{s=0}^{j-1}v_{j-s-1}v_{s}, j=1,2,\ldots, \; v_{0}=\left\|u_{n}^{(0)}(x)\right\|=1.
\end{equation}
It is easy to see that the sequence $\left\{V_{i}\right\},$ defined by the following recurrence formula,
\begin{equation}\label{gen_alg_24}
    V_{j}=\sum\limits_{s=0}^{j-1}V_{j-1-s}V_{s},\quad V_{0}=1,\quad j=0,1,\ldots,
\end{equation}
satisfies inequalities
\begin{equation}\label{gen_alg_25}
    V_{j}\geq v_{j},\quad j=0,1,\ldots.
\end{equation}
Let us consider the generating function $f(z)$ defined by the formula
\begin{equation}\label{gen_alg_26}
    f(z)=\sum\limits_{s=0}^{\infty}V_{s}z^{s}, \quad z\in \mathbb{R}.
\end{equation}
From \eqref{gen_alg_24} we obtain the equation with respect to $f(z)$
\begin{equation}\label{gen_alg_27}
    f(z)=z f^{2}(z)+1,\quad z\in (-\rho, \rho),
\end{equation}
where $\rho$ denotes the convergence radius for the power series \eqref{gen_alg_26}.
The solution to equation \eqref{gen_alg_27}, which satisfies the condition $f(0)=1$ is
\begin{equation}\label{gen_alg_28}
    f(z)=(2z)^{-1}\left(1-\sqrt{1-4z}\right)=\left(2z+\sum\limits_{s=2}^{\infty}\frac{(2s-3)!!}{(2s)!!}(4z)^{s}\right)(2z)^{-1},\quad z\in \left[-\frac{1}{4}, \frac{1}{4}\right],
\end{equation}
and
\begin{equation}\label{gen_alg_29}
    V_{j}=2\frac{(2j-1)!!}{(2j+2)!!}4^{j}=\alpha_{j}4^{j},\quad j=1,2,\ldots, \quad (-1)!!\stackrel{\textmd{def}}{=}1.
\end{equation}
\eqref{gen_alg_22}, \eqref{gen_alg_25} together with \eqref{gen_alg_29} lead us to the estimation
\begin{equation}\label{gen_alg_30}
  \left\|u_{n}^{(j)}(x)\right\|=\left(M_{n}\left\|q(x)-\bar{q}(x)\right\|_{\infty, [-1,1]}\right)^{j}v_{j}\leq \left(4 \bar{r}_{n}\right)^{j}\alpha_{j}=r_{n}^{j}\alpha_{j}.
\end{equation}
And inequality \eqref{gen_alg_30} yields the following error estimations of the FD-method's convergence rate
\begin{equation}\label{gen_alg_31}
    \left\|u_{n}(x)-\stackrel{m}{u}_{n}(x)\right\|\leq
    \left\{
    \begin{array}{cc}
    \cfrac{r_{n}^{m+1}}{1-r_{n}}\alpha_{m+1}, & \mbox{if} \quad r_{n}<1 \\
    \sum\limits_{j=m+1}^{\infty}\cfrac{1}{(j+1)\sqrt{\pi j}}, & \mbox{if}\quad r_{n}=1, \\
    \end{array}
    \right.
\end{equation}
\begin{equation}\label{gen_alg_32}
    \left|\lambda_{n}-\stackrel{m}{\lambda}_{n}\right|\leq \left\{
    \begin{array}{cc}
    \left\|q(x)-\bar{q}(x)\right\|_{\infty, [-1,1]}\cfrac{r_{n}^{m}}{1-r_{n}}\alpha_{m}, & \mbox{if}\quad r_{n}<1, \\
    \left\|q(x)-\bar{q}(x)\right\|_{\infty, [-1,1]}\sum\limits_{j=m}^{\infty}\cfrac{1}{(j+1)\sqrt{\pi j}}, & \mbox{if}\quad r_{n}=1. \\
    \end{array}
    \right.
\end{equation}
The infinite sum in the right sides of inequalities \eqref{gen_alg_31} and \eqref{gen_alg_32} can be easily estimated
\begin{equation}\label{gen_alg_31_est}
    \sum\limits_{j=m}^{\infty}\cfrac{1}{(j+1)\sqrt{\pi j}}\leq \int\limits_{m-1}^{\infty}\frac{1}{(x+1)\sqrt{x\pi}} d x=\frac{2}{\sqrt{\pi}}\arctan\left(\frac{1}{\sqrt{m+1}}\right).
\end{equation}

We have proved the following theorem.
\begin{theorem}\label{main_theorem}
   Let $q(x) \in Q^{0}[-1,1]$ and the constant $\bar{r}_{n}=M_{n}\left\|q(x)-\bar{q}(x)\right\|_{\infty, [-1,1]}$ satisfies the inequality $\bar{r}_{n}\leq\frac{1}{4}.$ Then the FD-method \eqref{gen_alg_3}, \eqref{gen_alg_5} -- \eqref{gen_alg_7_mk} converges faster than the geometric series with denominator $\bar{r}_{n}$ (superexponentially). The estimations of its convergence rate are given by formulas \eqref{gen_alg_31}, \eqref{gen_alg_32}.
\end{theorem}

Before passing to the next section, let us investigate the question about the asymptotical  behavior of parameter $M_{n}$ \eqref{gen_alg_20_Mn} as $n$ tends to $+\infty.$ It is easy to see that the basic problem \eqref{gen_alg_8} -- \eqref{gen_alg_9_mk_bp} is equivalent to  the following auxiliary eigenvalue problem
\begin{equation}\label{gen_alg_8_aux}
    \frac{\partial}{\partial x}\left[(1-x^{2})\frac{\partial v(x, \tau)}{\partial x}\right]+\left[\mu(\tau)-\tau\bar{q}(x)\right]v(x, \tau)=0,\quad x\in (-1,1),\; \tau\in [0,1],
\end{equation}
\begin{equation}\label{gen_alg_9_aux}
    \lim\limits_{x\rightarrow \pm 1}\left[\left(1-x^{2}\right)\frac{\partial v(x,\tau)}{\partial x}\right]=0,\quad \int\limits_{-1}^{1}\left(v(x, \tau)\right)^{2} d x=1,\; \forall \tau\in [0,1],
\end{equation}
\begin{equation}\label{gen_alg_9_mk_bp_aux}
    \left[v(x,\tau)\right]_{x=x_{i}}=0,\quad \left[\frac{\partial v(x, \tau)}{\partial x}\right]_{x=x_{i}}=0,\;\forall \tau\in [0,1],\; i=1,2,\ldots, N-1,
\end{equation}
when $\tau=1.$ On the other hand, as it follows from the proof of Theorem \ref{teor_1} (see \cite{Klimenko_Ya}), for sufficiently large number $n,$ namely
\begin{equation}\label{gen_alg_33_n_estimation}
    n>2\left\|q(x)\right\|_{\infty, [-1,1]}\geq 2\left\|\bar{q}(x)\right\|_{\infty, [-1,1]},
\end{equation}
the eigensolution $\mu_{n}(\tau), \; v_{n}(x, \tau)$ to the problem \eqref{gen_alg_8_aux} -- \eqref{gen_alg_9_mk_bp_aux} exists and can be developed to the power series with respect to $\tau$ on the interval $[0, 1],$ for all $x\in (-1,1),$ furthermore
\begin{equation}\label{gen_alg_33_mu_cond}
   \mu_{n}(0)=n(n+1), \mu_{n}(1)=\lambda_{n}^{(0)}.
\end{equation}
 Hence, from equality \eqref{gen_alg_8_aux} it is easy to obtain
\begin{equation}\label{gen_alg_33}
    \frac{\partial}{\partial x}\left[(1-x^{2})\frac{\partial^{2} v_{n}(x, \tau)}{\partial x\partial \tau}\right]+\left[\mu_{n}(\tau)-\tau\bar{q}(x)\right]\frac{\partial}{\partial \tau}v_{n}(x, \tau)=\left[\mu^{\prime}_{n}(\tau)-\bar{q}(x)\right]v_{n}(x, \tau),
\end{equation}
\noindent $\quad x\in (-1,1),\; \tau\in [0,1].$ After that, applying to the both sides of equality \eqref{gen_alg_33} the integral operator $\int\limits_{-1}^{1}v_{n}(x, \tau)(\cdot)d x,$ we obtain
\begin{equation}\label{gen_alg_34}
    \mu^{\prime}_{n}(\tau)=\int\limits_{-1}^{1}\left(v_{n}(x, \tau)\right)^{2}\bar{q}(x)d x,\quad \tau\in [0,1].
\end{equation}
Taking into account \eqref{gen_alg_33_mu_cond} and equality \eqref{gen_alg_34}, it is easy to obtain
\begin{equation}\label{gen_alg_35}
    \lambda_{n}^{(0)}=n(n+1)+\int\limits_{0}^{1}\int\limits_{-1}^{1}\left(v_{n}(x, \tau)\right)^{2}\bar{q}(x)d x d \tau.
\end{equation}
Finally, equality \eqref{gen_alg_35} yields the following estimation
\begin{equation}\label{gen_alg_36}
    \lambda_{n+1}^{(0)}-\lambda_{n}^{(0)}\geq 2(n+1)-2\left\|q(x)\right\|_{\infty, [-1,1]}>0
\end{equation}
Hence, for the $n$ sufficiently large (see \eqref{gen_alg_33_n_estimation}) we have the following estimation which describes the asymptotical behavior of the $M_{n}:$
\begin{equation}\label{gen_alg_37}
    M_{n}\leq \left(2n-2\left\|q(x)\right\|_{\infty, [-1,1]}\right)^{-1}.
\end{equation}

\section{General algorithm of the FD-method: software implementation.}

The solution to the basic problem \eqref{gen_alg_8} -- \eqref{gen_alg_9_mk_bp} can be represented in the following form
\begin{equation}\label{alg_implement_1}
  u_{n}^{(0)}(x)=A_{i}P_{\nu_{i}}(x)+B_{i}Q_{\nu_{i}}(x),\quad x\in \left[x_{i-1}, x_{i}\right],\quad A_{i}, B_{i}\in \mathbb{R},\quad i=\overline{1, N},
\end{equation}
where
\begin{equation}\label{alg_implement_1_nu}
    \nu_{i}=\nu_{i}\left(\lambda_{n}^{(0)}\right)=\frac{1}{2}\left(-1\pm\sqrt{1+4\left(\lambda^{(0)}_{n}-\bar{q}_{i}\right)}\right)\footnote{The sing ``$+$'' or ``$-$'' can be chosen optionally.},\quad \bar{q}_{i}\equiv \bar{q}(x),\; x\in \left[x_{i-1}, x_{i}\right)
\end{equation}
and $Q_{\nu_{i}}(x)$ denotes the Legendre function of the second kind (see \cite{Bateman_1}).
It is well known that (see, for example, \cite{Bateman_1})
\begin{eqnarray}\label{alg_implement_2}
& \lim\limits_{x\rightarrow -1}\left(1-x^{2}\right)\cfrac{d P_{\nu}(x)}{d x}=\cfrac{2\sin\left(\pi \nu\right)}{\pi}, \quad \lim\limits_{x\rightarrow 1}\left(1-x^{2}\right)\cfrac{d P_{\nu}(x)}{d x}=0, \nonumber\\
\\
& \lim\limits_{x\rightarrow -1}\left(1-x^{2}\right)\cfrac{d Q_{\nu}(x)}{d x}=\cos\left(\pi \nu\right), \quad \lim\limits_{x\rightarrow 1}\left(1-x^{2}\right)\cfrac{d Q_{\nu}(x)}{d x}=1.\nonumber
\end{eqnarray}

To satisfy boundary condition \eqref{gen_alg_9} we have to set
\begin{equation}\label{alg_implement_3}
 A_{1}\cfrac{2\sin\left(\pi \nu_{1}\right)}{\pi}+B_{1}\cos\left(\pi \nu_{1}\right)=0, \quad B_{N}=0.
\end{equation}
The other constants $A_{i}, B_{i}$ have to be found from the matching conditions \eqref{gen_alg_9_mk_bp} as the solutions to the following recurrence systems of linear algebraic equations:
\begin{equation}\label{alg_implement_4}
    \left\{
      \begin{array}{ccc}
        A_{i-1}P_{\nu_{i-1}}(x_{i-1})+B_{i-1}Q_{\nu_{i-1}}(x_{i-1}) &=& A_{i}P_{\nu_{i}}(x_{i-1})+B_{i}Q_{\nu_{i}}(x_{i-1}), \\
        A_{i-1}P^{\prime}_{\nu_{i-1}}(x_{i-1})+B_{i-1}Q^{\prime}_{\nu_{i-1}}(x_{i-1}) &=& A_{i}P^{\prime}_{\nu_{i}}(x_{i-1})+B_{i}Q^{\prime}_{\nu_{i}}(x_{i-1}), \\
      \end{array}
    \right.
\end{equation}
$i=2,3,\ldots, N.$

Let us consider the function $$\Phi\left(\lambda\right)=A_{1}\cfrac{2\sin\left(\pi \nu_{1}\right)}{\pi}+B_{1}\cos\left(\pi \nu_{1}\right),$$ where constants $A_{1}, B_{1}$ are defined by the following recurrence formulas (see \eqref{alg_implement_4})
\begin{eqnarray}\label{alg_implement_5}
   &A_{i-1}= \cfrac{1}{\delta_{i-1}}\left(C_{i}Q^{\prime}_{\nu_{i-1}}(x_{i-1})-D_{i}Q_{\nu_{i-1}}(x_{i-1})\right), \nonumber\\
   &B_{i-1}= \cfrac{1}{\delta_{i-1}}\left(D_{i}P_{\nu_{i-1}}(x_{i-1})-C_{i}P^{\prime}_{\nu_{i-1}}(x_{i-1})\right), \nonumber\\
   &\delta_{i-1}= P_{\nu_{i-1}}(x_{i-1})Q^{\prime}_{\nu_{i-1}}(x_{i-1})-P^{\prime}_{\nu_{i-1}}(x_{i-1})Q_{\nu_{i-1}}(x_{i-1}), \\
   &C_{i}= A_{i}P_{\nu_{i}}(x_{i-1})+B_{i}Q_{\nu_{i}}(x_{i-1}), \nonumber\\
   &D_{i}= A_{i}P^{\prime}_{\nu_{i}}(x_{i-1})+B_{i}Q^{\prime}_{\nu_{i}}(x_{i-1}),\quad i=2,3,\ldots, N \nonumber
\end{eqnarray}
with the initial data
\begin{equation}\label{alg_implement_5_in_data}
    A_{N}=1,\quad B_{N}=0,
\end{equation}
 and the values $\nu_{i}=\nu_{i}\left(\lambda\right)$ are defined by formula \eqref{alg_implement_1_nu}.
It is easy to see that the eigenvalues of the basic problem \eqref{gen_alg_8} -- \eqref{gen_alg_9_mk_bp} coincide with the solutions of the equation
\begin{equation}\label{alg_implement_6}
    \Phi(\lambda)=0.
\end{equation}
Since the eigenvalues of the basic problem are simple, they can be found with any precision as the solutions to equation \eqref{alg_implement_6} by the bisection method.

Let $\lambda_{n}^{(0)}$ be an arbitrary solution to equation \eqref{alg_implement_6}. Using the coefficients $A_{i}, B_{i}, i=1,2,\ldots, N,$ computed by formulas \eqref{alg_implement_5} with $\lambda=\lambda_{n}^{(0)}$ we can compute corresponding eigenfunction $u_{n}^{(0)}(x)$ by formula \eqref{alg_implement_1}.

To compute $\lambda_{n}^{(j)}$ we can use the formula
\begin{equation}\label{alg_implement_7}
  \lambda_{n}^{(j)}=\left[\int\limits_{-1}^{1}\left(u_{n}^{(0)}(x)\right)^{2}d x\right]^{-1}\int\limits_{-1}^{1}\left(q(x)-\bar{q}(x)\right)u_{n}^{(j-1)}(x)u_{n}^{(0)}(x)d x.
\end{equation}
 And for the computation of the functions $u_{n}^{(j)}(x), \; j=1,2,\ldots, m$ it is convenient to use the formula
\begin{equation}\label{alg_implement_7_1}
   u_{n}^{(j)}(x)=E_{j}u_{n}^{(0)}(x)+\int\limits_{-1}^{x}K(x,\xi)F^{(j)}_{n}(\xi) d \xi=
\end{equation}
$$=E_{j}u_{n}^{(0)}(x)+w_{n}(x)\int\limits_{-1}^{x}u_{n}^{(0)}(\xi)F^{(j)}_{n}(\xi) d \xi-u_{n}^{(0)}(x)\int\limits_{-1}^{x}w_{n}(\xi)F^{(j)}_{n}(\xi) d \xi,$$
where
$$K(x,\xi)=w_{n}(x)u_{n}^{(0)}(\xi)-w_{n}(\xi)u_{n}^{(0)}(x),$$
$$E_{j}=-\left[\int\limits_{-1}^{1}\left(u_{n}^{(0)}(x)\right)^{2}d x\right]^{-1}\int\limits_{-1}^{1}\left\{u_{n}^{(0)}(x)\int\limits_{-1}^{x}K(x,\xi)F^{(j)}_{n}(\xi) d \xi\right\} d x$$
and
$w_{n}(x)$ denotes a function, defined by the formula
$$  w_{n}(x)=A_{i}P_{\nu_{i}}(x)+B_{i}Q_{\nu_{i}}(x),\quad x\in \left[x_{i-1}, x_{i}\right], \quad i=\overline{1, N},$$
where $A_{i}, B_{i},\; i=1,2,\ldots, N$ can be computed recursively by formulas \eqref{alg_implement_5} with the initial data
\begin{equation}\label{alg_implement_8}
    A_{N}=0,\;B_{N}=1.
\end{equation}
Taking into account the properties of the Legendre functions (see \cite{Bateman_1}), it is easy to verify that
\begin{equation}\label{alg_implement_9}
  \left.\frac{\partial K(x,\xi)}{\partial x}\right|_{\xi=x}=\left(1-x^{2}\right)^{-1},\quad x\in(-1,1)
\end{equation}
and $K(x, \xi)$ is the Cauchy operator for the nonhomogeneous equations \eqref{gen_alg_5}.

In general case, the integrals in formulas \eqref{alg_implement_7} and \eqref{alg_implement_7_1} could not be calculated analytically. Moreover, the integrands in this formulas are unbounded at the points $\pm 1.$ To compute this integrals it is convenient to use numerical methods, for example, the {\it tanh} rule (see \cite{Kahaner_Nesh}) and Stenger's formula (see \cite{Stenger_num_meth}):
\begin{equation}\label{alg_implement_10}
    \int_{a}^{b}f(x)d x=\int\limits_{-\infty}^{+\infty}f\left(\frac{a+be^{\omega}}{1+e^{\omega}}\right)\frac{\left(b-a\right)d\omega}{\left(e^{-\omega/2}+e^{\omega/2}\right)^{2}}\approx
\end{equation}
$$\approx h_{sinc}\sum\limits_{l=-K}^{K}f\left(\frac{a+be^{lh_{sinc}}}{1+e^{lh_{sinc}}}\right)\frac{\left(b-a\right)}{\left(e^{-lh_{sinc}/2}+e^{lh_{sinc}/2}\right)^{2}}$$
and
\begin{equation}\label{alg_implement_11}
   \int_{a}^{z_{k}}f(x)d x\approx h_{sinc}\sum\limits_{l=-K}^{K}\delta_{k-l}^{(-1)}f\left(\frac{a+be^{lh_{sinc}}}{1+e^{lh_{sinc}}}\right)\frac{\left(b-a\right)}{\left(e^{-lh_{sinc}/2}+e^{lh_{sinc}/2}\right)^{2}}\footnote{ The function $f(x)$ is required to be sufficiently smooth on (a, b), see \cite{Stenger_num_meth}. },
\end{equation}
where
$$z_{k}=\frac{a+be^{h_{sinc}k}}{1+e^{h_{sinc}k}},\;k=-K \ldots, K,\quad \delta_{k}^{(-1)}=\frac{1}{2}+\int\limits_{0}^{k}\frac{\sin\left(\pi t\right)}{\pi t}d t, \; k=-2K\ldots 2K,$$
and $h_{sinc}=\sqrt{\frac{2\pi}{K}}\footnote{About the optimal choice of the value for $h_{sinc}$ see \cite{Stenger_num_meth}.}$

Henceforth we require that the function $q(x)$ is analytical on each interval $(x_{i-1}, x_{i}),$ $i=1,2,\ldots, N.$
Before passing to the algorithm description, we need to introduce the following auxiliary notation
\begin{equation}\label{alg_implement_12}
    z_{i,j}=\frac{x_{i-1}+x_{i}e^{h_{sinc}j}}{1+e^{h_{sinc}j}},\quad \mu_{i,j}=\frac{\left(x_{i-1}-x_{i}\right)}{\left(e^{-jh_{sinc}/2}+e^{jh_{sinc}/2}\right)^{2}},\; \quad \nu_{n}=\left[\int\limits_{-1}^{1}\left(u_{n}^{(0)}(x)\right)^{2}d x\right]^{-1},
\end{equation}
$$i=1,2,\ldots, N,\;j=-K,\ldots, K.$$
One of the possible algorithms for computation of $\lambda_{n}^{(k)}$ and $u_{n}^{(k)}(k)$ for $k=1,2,\ldots, r,$ $i=1,2,\ldots, N,$ $j=-K\ldots, K$ is described below.
\\

\RestyleAlgo{ruled}
\DontPrintSemicolon
\begin{algorithm}[H]
\RestyleAlgo{ruled}
\SetKwInOut{Input}{input}\SetKwInOut{Output}{output}
  \Input{The arrays of values $u_{n}^{(0)}(z_{i,j}),\; w_{n}(z_{i, j}),\;\mu_{i,j},\quad i\in \overline{1,N},\; j\in\overline{-K,K},$ $\delta_{l}^{(-1)},\; l\in\overline{-2K..2K},$
  $\lambda_{n}^{0},$ and $r$ --- the depth of the FD-method.}
  \Output{The arrays of values $u_{n}^{(d)}(z_{i,j}),\; d\in\overline{1, r},\;i\in \overline{1,N},\quad j\in\overline{-K,K},$
  $\lambda_{n}^{(d)},\; d\in\overline{1,r}.$
  }

\Begin{
\nl  \For{$d:=1$ \KwTo $r$}{
\nl  {\bf ComputeLambda( $u_{n}^{(d)}{(z_{i,j})}$, $u_{n}^{(0)}{(z_{i,j})}$)}\tcp{This procedure computes $\lambda_{n}^{(d)},$ see description below}
\nl {\bf ComputeF($u_{n}^{(k)}{(z_{i,j})},$ $\lambda_{n}^{(k)}$)}\tcp{This procedure computes the values $F_{n}^{(d)}(z_{i,j})$, see description below}
\nl {\bf ComputeCorrectionForEigenfunction($u_{n}^{(0)}{(z_{i,j})},$ $w_{n}{(z_{i,j})},$ $F_{n}^{(d)}{(z_{i,j})}$)}\;\tcp{This procedure computes the values $u_{n}^{(d)}(z_{i,j}),$  see description below}
\nl {\bf ProvideOrthogonalityCondition( $u_{n}^{(d)}{(z_{i,j})}$, $u_{n}^{(0)}{(z_{i,j})}$)}\tcp{This procedure recompute the values $u_{n}^{(d)}{(z_{i,j})}$ to provide the orthogonality condition $\int\limits_{-1}^{1}u_{n}^{(0)}(x)u_{n}^{(d)}(x)d x=0$}
  }
  }
%\SetAlgoProcName{aname}
\caption{The algorithm for computation of $\lambda_{n}^{d},$ $u_{n}^{d}(z_{i,j}),$ $d\in\overline{1, r},\;i\in\overline{1, N},\; j\in \overline{-K, K}$}
\end{algorithm}

\begin{procedure}
\caption{ComputeLambda( $u_{n}^{(d)}{(z_{i,j})}$, $u_{n}^{(0)}{(z_{i,j})}$)}
\SetKwInOut{Input}{input}\SetKwInOut{Output}{output}
\Input{The arrays of values $u_{n}^{(d)}{(z_{i,j})}$, $u_{n}^{(0)}{(z_{i,j})},$ \quad $i\in\overline{1, N}, \; j\in \overline{-K,K}$}
 \Output{ The value $\lambda_{n}^{(d)}$}
\Begin{
\nl $\lambda_{n}^{(d)}:=0$\;
\nl     \For{$i:=1$ \KwTo $N$}{
\nl        \For{$j:=-K$ \KwTo $K$}{
\nl          $\lambda_{n}^{(d)}:=\lambda_{n}^{(d)}+u_{n}^{(0)}(z_{i,j})u_{n}^{(d-1)}(z_{i,j})\left\{q(z_{i,j})-\bar{q}(z_{i,j})\right\}\mu_{i,j}$ \tcp{Compute $\lambda_{n}^{(d)}$ by formula \eqref{alg_implement_7}, using formula \eqref{alg_implement_10}}
        }
     }
\nl     $\lambda_{n}^{(d)}:=h_{sinc}\nu_{n}\lambda_{n}^{(d)}$\tcp{The $\lambda_{n}^{(d)}$ is computed}
}
\end{procedure}

\begin{procedure}[H]
\caption{ProvideOrthogonalityCondition( $u_{n}^{(d)}{(z_{i,j})}$, $u_{n}^{(0)}{(z_{i,j})}$)}\label{Lambda_comp}
\SetKwInOut{Input}{input}\SetKwInOut{Output}{output}
\Input{The arrays of values $u_{n}^{(d)}{(z_{i,j})}$, $u_{n}^{(0)}{(z_{i,j})},$ \quad $i\in\overline{1, N}, \; j\in \overline{-K,K}$}
 \Output{The array of values $u_{n}^{(d)}{(z_{i,j})},$ \quad $i\in\overline{1, N}, \; j\in \overline{-K,K},$ such that the integral $\int\limits_{-1}^{1}u_{n}^{(0)}(x)u_{n}^{(d)}(x),$ computed by formula \eqref{alg_implement_10}, is equal to zero}
\Begin{
\nl$\mathcal{I}:=0$\;
\nl  \For{$i:=1$ \KwTo $N$}{
\nl    \For{$j:=-K$ \KwTo $K$}{
\nl      $\mathcal{I}:=\mathcal{I}+u_{n}^{(d)}(z_{i,j})u_{n}^{(0)}(z_{i,j})\mu_{i,j}$
    }
  }
\nl  $\mathcal{I}:=\mathcal{I}h_{sinc}\nu_{n}$\;
\nl  \For{$i=1$ \KwTo $N$}{
\nl    \For{$j=-K$ \KwTo $K$}{
\nl      $u_{n}^{(d)}(z_{i,j}):=u_{n}^{(d)}(z_{i,j})-\mathcal{I}u_{n}^{(0)}(z_{i,j})$\;
    }
  }
}
\end{procedure}

\begin{procedure}
\caption{ComputeF($u_{n}^{(k)}{(z_{i,j})},$ $\lambda_{n}^{(k)}$)}
\SetKwInOut{Input}{input}\SetKwInOut{Output}{output}
\Input{The arrays of values $u_{n}^{(k)}{(z_{i,j})}$, \quad $k\in\overline{0, d},\;i\in\overline{1, N}, \; j\in \overline{-K,K}$ and $\lambda_{n}^{(k)},$ $k\in\overline{0,d}$}
\Output{The array of values $F_{n}^{(d)}(z_{i,j}),$ $i\in\overline{1, N}, \; j\in \overline{-K,K},$ computed by formula \eqref{gen_alg_5}}
\Begin{
\nl\For{$i:=1$ \KwTo $N$}{
\nl  \For{$j:=-K$ \KwTo $K$}{
\nl       $F_{n}^{(d)}(z_{i,j}):=0$\;
\nl       \For{$k:=0$ \KwTo $d-1$}{
\nl          $F_{n}^{(d)}(z_{i,j}):=F_{n}^{(d)}(z_{i,j})-\lambda_{d-k}u_{n}^{(k)}(z_{i,j})$
       }
\nl       $F_{n}^{(d)}(z_{i,j}):=F_{n}^{(d)}(z_{i,j})+u_{n}^{(d-1)}(z_{i,j})\left\{q(z_{i,j})-\bar{q}(z_{i,j})\right\}$
  }
}

}
\end{procedure}

\begin{procedure}[H]
\caption{ComputeCorrectionForEigenfunction($u_{n}^{(0)}{(z_{i,j})},$ $w_{n}{(z_{i,j})},$ $F_{n}^{(d)}{(z_{i,j})}$)}

\SetKwInOut{Input}{input}\SetKwInOut{Output}{output}
\Input{The arrays of values $u_{n}^{(0)}{(z_{i,j})}$, $w_{n}{(z_{i,j})}$ and $F_{n}^{(d)}{(z_{i,j})},$ \quad $i\in\overline{1, N}, \; j\in \overline{-K,K}$}
 \Output{The array of values $u_{n}^{(d)}{(z_{i,j})},$ \quad $i\in\overline{1, N}, \; j\in \overline{-K,K}$}

\Begin{
\nl  $\mathcal{V}_1:=0;$ $\mathcal{V}_2:=0$\tcp{Initialize auxiliary variables}
\nl  \For{$i:=1$ \KwTo $N$}{
\nl    \For{$j:=-K$ \KwTo $K$}{
\nl      $\mathcal{I}_{1}:=0;$ $\mathcal{I}_{2}:=0$\;
\nl      \For{$l:=-K$ \KwTo $K$}{
\nl      $\mathcal{I}_{1}:=\mathcal{I}_{1}+u_{n}^{(0)}(z_{i,l})\delta_{j-l}^{(-1)}F^{(d)}(z_{i,l})\mu_{i,l}$\tcp{Compute integral $\int\limits_{-1}^{z_{i,j}}u_{n}^{(0)}(\xi)F^{(j)}_{n}(\xi) d \xi$ by the Stenger's formula \eqref{alg_implement_11}}
\nl      $\mathcal{I}_{2}:=\mathcal{I}_{2}+w_{n}(z_{i,l})\delta_{j-l}^{(-1)}F^{(d)}(z_{i,l})\mu_{i,l}$\tcp{Compute integral $\int\limits_{-1}^{z_{i,j}}w_{n}(\xi)F^{(j)}_{n}(\xi) d \xi$ by the Stenger's formula \eqref{alg_implement_11}}
      }
\nl      $\mathcal{I}_{1}:=\mathcal{I}_{1}h_{sinc}+\mathcal{V}_{1}$\;
\nl      $\mathcal{I}_{2}:=\mathcal{I}_{2}h_{sinc}+\mathcal{V}_{2}$\;
\nl      $u_{n}^{(d)}(z_{i,j}):=w_{n}(z_{i,j})\mathcal{I}_{1}-u_{n}^{(0)}(z_{i,j})\mathcal{I}_{2}$\;
    }
\nl     $\mathcal{V}_{1}:=\mathcal{I}_{1}$\;
\nl     $\mathcal{V}_{2}:=\mathcal{I}_{2}$\;
  }

}
\end{procedure}

The described algorithm does not clime to be the most efficient. It is evident that the highest accuracy, that can be achieved by increasing the rank $r$ of the FD-method, is restricted by the accuracy of the quadrature formulas \eqref{alg_implement_10}, \eqref{alg_implement_11}, used in this algorithm. The question about optimal choice of the parameters $r$ and $K$ is still a pressing issue. In the next section the numerical results obtained by the Algorithm 1 are stated. In the first example we apply quadrature formulas \eqref{alg_implement_10}, \eqref{alg_implement_11} with $K=500$ and in the other examples we used $K=350.$

\section{Numerical examples}

{\bf Example 1.} Let us consider the following Sturm-Liouville problem
\begin{equation}\label{examp_1}
    L\left(u(x)\right)=\frac{d}{d x}\left[\left(1-x^{2}\right)\frac{d u(x)}{d x}\right]+\left(\lambda-q(x)\right)u(x)=0,\quad x\in (-1,1),
\end{equation}
$$\lim\limits_{x\rightarrow \pm 1}\left[\left(1-x^{2}\right)\frac{d u(x)}{d x}\right]=0$$
with
\begin{equation}\label{examp_5}
q(x)=x.
\end{equation}

Using the FD-method, described in \cite{Klimenko_Ya} and the computer algebra system Maple, it is easy to obtain
$$\left.
  \begin{array}{ll}
\lambda_{0}^{(0)}=0,&u_{0}^{(0)}(x)=\cfrac{\sqrt{2}}{2},\\
   \lambda_{0}^{(1)}=0,&u_{0}^{(1)}(x)=-\cfrac{\sqrt{2}x}{4},\\
   \lambda_{0}^{(2)}=-\cfrac{1}{6},&u_{0}^{(2)}(x)=\cfrac{\sqrt{2}x^{2}}{24}-\cfrac{\sqrt{2}}{72},\\
   \lambda_{0}^{(3)}=0,&u_{0}^{(3)}(x)=-\cfrac{\sqrt{2}x^{3}}{288}+\cfrac{5\sqrt{2} x}{288},\\
   \lambda_{0}^{(4)}=\cfrac{11}{1080},&u_{0}^{(4)}(x)=\cfrac{\sqrt{2} x^{4}}{5760}-\cfrac{\sqrt{2}x^{2}}{270}+\cfrac{311\sqrt{2}}{259200},\\
   \lambda_{0}^{(5)}=0,&u_{0}^{(5)}(x)=-\cfrac{\sqrt{2}x^{5}}{172800}+\cfrac{\sqrt{2}x^{3}}{2880}-\cfrac{ 1181\sqrt{2} x}{518400},\\
   \lambda_{0}^{(6)}=-\cfrac{47}{34020},&u_{0}^{(6)}(x)=\cfrac{\sqrt{2}x^{6}}{7257600}-\cfrac{\sqrt{2}x^{4}}{53760}+\cfrac{11237\sqrt{2}x^{2}}{21772800}-\cfrac{76967\sqrt{2}}{457228800}.  \end{array}
\right.
$$

To control the accuracy of the results obtained by the FD-method we use the following functionals
\begin{equation}\label{examp_2}
    \stackrel{m}{\eta}_{n}=\left[\int\limits_{-1}^{1}\left[\left(1-x^{2}\right)\frac{d \stackrel{m}{u}_{n}(x)}{d x}+\int\limits_{-1}^{x}\left(\stackrel{m}{\lambda}_{n}-q(\xi)\right)\stackrel{m}{u}_{n}(\xi)d\xi\right]^{2} d x\right]^{\frac{1}{2}},
\end{equation}

\begin{equation}\label{examp_2}
   \stackrel{m}{\bar{\eta}}_{n}=\left[\int\limits_{-1}^{1}\left[\frac{d}{d x}\left[\left(1-x^{2}\right)\frac{d \stackrel {m}{u}_{n}(x)}{d x}\right]+\left(\stackrel {m}{\lambda}_{n}-q(x)\right)\stackrel {m}{u}_{n}(x)\right]^{2} d x\right]^{\frac{1}{2}}
\end{equation}

Also we compare the FD-method results with eigenvalues obtained by SLAIGN2 (see. \cite{Alg_810})

\begin{table}[htbp]
\caption{Example 1, the data obtained by SLEIGN2}
\centering
\begin{tabular}{|c|c|c|c|}
  \hline
  $n$ & $\lambda_{n, sl2}$ & $TOL$ & $IFLAG$ \\
  \hline
  0 & $-0.157663483D+00$ & $0.23950D-13$ & $1$ \\
  \hline
  1 & $0.209076065D+01$ & $0.19238D-13$ & $1$ \\
  \hline
  2 & $0.602403235D+01$ & $0.20718D-08$ & $1$ \\
  \hline
  3 & $0.120111226D+02$ & $0.10472D-12$ & $1$ \\
  \hline
  4 & $0.200064954D+02$ & $0.70433D-13$ & $1$ \\
  \hline
\end{tabular}
\label{table_1}
\end{table}
In all examples, presented in this section, we apply the FD-method with uniform mesh on the interval $[-1, 1]$ and $N$ denotes the number of subintervals.
The results obtained by the FD-method are presented in tables \ref{table_2}-- \ref{table_4} and figure \ref{image1} below. The result obtained by the SLEIGN2 are presented in table \ref{table_1}.

\begin{table}[htbp]
\caption{Example 1, the $0$-th eigenvalue, FD-method with $N=1$}
\centering
\begin{tabular}{|c|c|c|c|c|}
  \hline
  $m$ & $\stackrel{m}{\lambda}_{0}$ & $\left\|u_{0}^{(m)}(x)\right\|$ & $\stackrel{m}{\bar{\eta}}_{n}$ & $\left|\stackrel {m}{\lambda}_{0}-\lambda_{0, sl2}\right|$ \\
  \hline
  0 & $0.0$ & $1.0$ & $0.56$ & $0.157663483$ \\
  \hline
  1 & $0.0$ & $0.29$ & $0.22$ & $0.157663483$ \\
  \hline
  2 & $-0.1666666667$ & $0.24e-1$ & $0.36\times 10^{-1}$ & $0.90031875\times 10^{-2}$ \\
  \hline
  3 & $-0.1666666667$ & $0.18\times 10^{-1}$ & $0.16\times 10^{-1}$ & $0.90031875\times 10^{-2}$ \\
  \hline
  4 & $-0.1564814815$ & $0.21\times 10^{-2}$ & $0.48\times 10^{-2}$ & $0.11819977\times 10^{-2}$ \\
  \hline
  5 & $-0.1564814815$ & $0.24\times 10^{-2}$ & $0.23\times 10^{-2}$ & $0.11819977\times 10^{-2}$ \\
  \hline
  6 & $-0.1578630218$ & $0.30\times 10^{-3}$ & $0.78\times 10^{-3}$ & $0.1995426\times 10^{-3}$ \\
  \hline
  7 & $-0.1578630218$ & $0.41\times 10^{-3}$ & $0.41\times 10^{-3}$ & $0.1995426\times 10^{-3}$ \\
  \hline
  8 & $-0.1576253633$ & $0.52\times 10^{-4}$ & $0.15\times 10^{-3}$ & $0.381159\times 10^{-4}$ \\
  \hline
  9 & $-0.1576253633$ & $0.79\times 10^{-4}$ & $0.78\times 10^{-4}$ & $0.381159\times 10^{-4}$ \\
  \hline
  10 & $-0.1576713252$ & $0.10\times 10^{-4}$ & $0.31\times 10^{-4}$ & $0.78460\times 10^{-5}$ \\
  \hline
\end{tabular}
\label{table_2}
\end{table}

\begin{table}[htbp]
\caption{Example 1, the $0$-th eigenvalue, FD-method with $N=1$}
\centering
\begin{tabular}{|c|c|c|c|c|}
  \hline
  $m$ & $\stackrel{m}{\lambda}_{0}$ & $\left\|u_{0}^{(m)}(x)\right\|$ & $\stackrel {m}{\bar{\eta}}_{n}$ & $\left|\stackrel {m}{\lambda}_{0}-\lambda_{0, sl2}\right|$ \\
  \hline
  51 & $-0.15766348313775096746$ & $0.12\times 10^{-16}$ & $0.13\times 10^{-16}$ & $0.39\times 10^{-8}$ \\
  \hline
  52 & $-0.15766348313775096031$ & $0.16\times 10^{-17}$ & $0.59\times 10^{-17}$ & $0.39\times 10^{-8}$ \\
  \hline
  53 & $-0.15766348313775096031$ & $0.33\times 10^{-17}$ & $0.33\times 10^{-17}$ & $0.39\times 10^{-8}$ \\
  \hline
  54 & $-0.15766348313775096218$ & $0.39\times 10^{-18}$ & $0.15\times 10^{-17}$ & $0.39\times 10^{-8}$ \\
  \hline
  55 & $-0.15766348313775096218$ & $0.84\times 10^{-18}$ & $0.86\times 10^{-18}$ & $0.39\times 10^{-8}$ \\
  \hline
  56 & $-0.15766348313775096169$ & $0.11\times 10^{-18}$ & $0.40\times 10^{-18}$ & $0.39\times 10^{-8}$ \\
  \hline
  57 & $-0.15766348313775096169$ & $0.22\times 10^{-18}$ & $0.22\times 10^{-18}$ & $0.39\times 10^{-8}$ \\
  \hline
  58 & $-0.15766348313775096182$ & $0.28\times 10^{-19}$ & $0.11\times 10^{-18}$ & $0.39\times 10^{-8}$ \\
  \hline
  59 & $-0.15766348313775096182$ & $0.58\times 10^{-19}$ & $0.61\times 10^{-19}$ & $0.39\times 10^{-8}$ \\
  \hline
  60 & $-0.15766348313775096178$ & $0.74\times 10^{-20}$ & $0.29\times 10^{-19}$ & $0.39\times 10^{-8}$ \\
  \hline
\end{tabular}
\label{table_3}
\end{table}

\begin{table}[htbp]
\caption{Example 1, FD-method with $N=3$}
\centering
\begin{tabular}{|c|c|c|c|c|c|c|}
  \hline
  $n$ & $m$ & $\stackrel{m}{\lambda}_{n}$ & $\left|\lambda_{n}^{(m)}\right|$ & $\left\|u_{n}^{(m)}(x)\right\|$ & $\stackrel {m}{\eta}_{n}$ & $\left|\stackrel{m}{\lambda}_{n}-\lambda_{n, sl2}\right|$ \\
  \hline
  0 & 15 & $-0.1576634831377509617898$ & $0.3\times 10^{-21}$ & $0.33\times 10^{-21}$ & $0.95\times 10^{-22}$ & $0.39\times 10^{-8}$ \\
  \hline
  1 & 15 & $2.090760648363956948786$ &$0.7\times 10^{-20}$ &$0.13\times 10^{-20}$ & $0.26\times 10^{-21}$ & $0.1482\times 10^{-5}$ \\
  \hline
  2 & 15 & $6.024031655336352711291$ & $0.7\times 10^{-20}$ &$0.94\times 10^{-21}$ & $0.22\times 10^{-21}$ & $0.5035\times 10^{-5}$ \\
  \hline
  3 & 13 & $12.01112256362987127625$ &$0.1\times 10^{-19}$ & $0.24\times 10^{-20}$ & $0.20\times 10^{-21}$ & $0.4\times 10^{-7}$ \\
  \hline
  4 & 11 & $20.00649533292656299628$ &$0.1\times 10^{-19}$ & $0.31\times 10^{-19}$ & $0.15\times 10^{-20}$ & $0.7\times 10^{-7}$ \\
  \hline
\end{tabular}
\label{table_4}
\end{table}

\begin{figure}[htbp]
\begin{minipage}[h]{1\linewidth}
\begin{minipage}[h]{0.48\linewidth}
\center{\rotatebox{-0}{\includegraphics[height=0.8\linewidth,
width=1.0\linewidth]{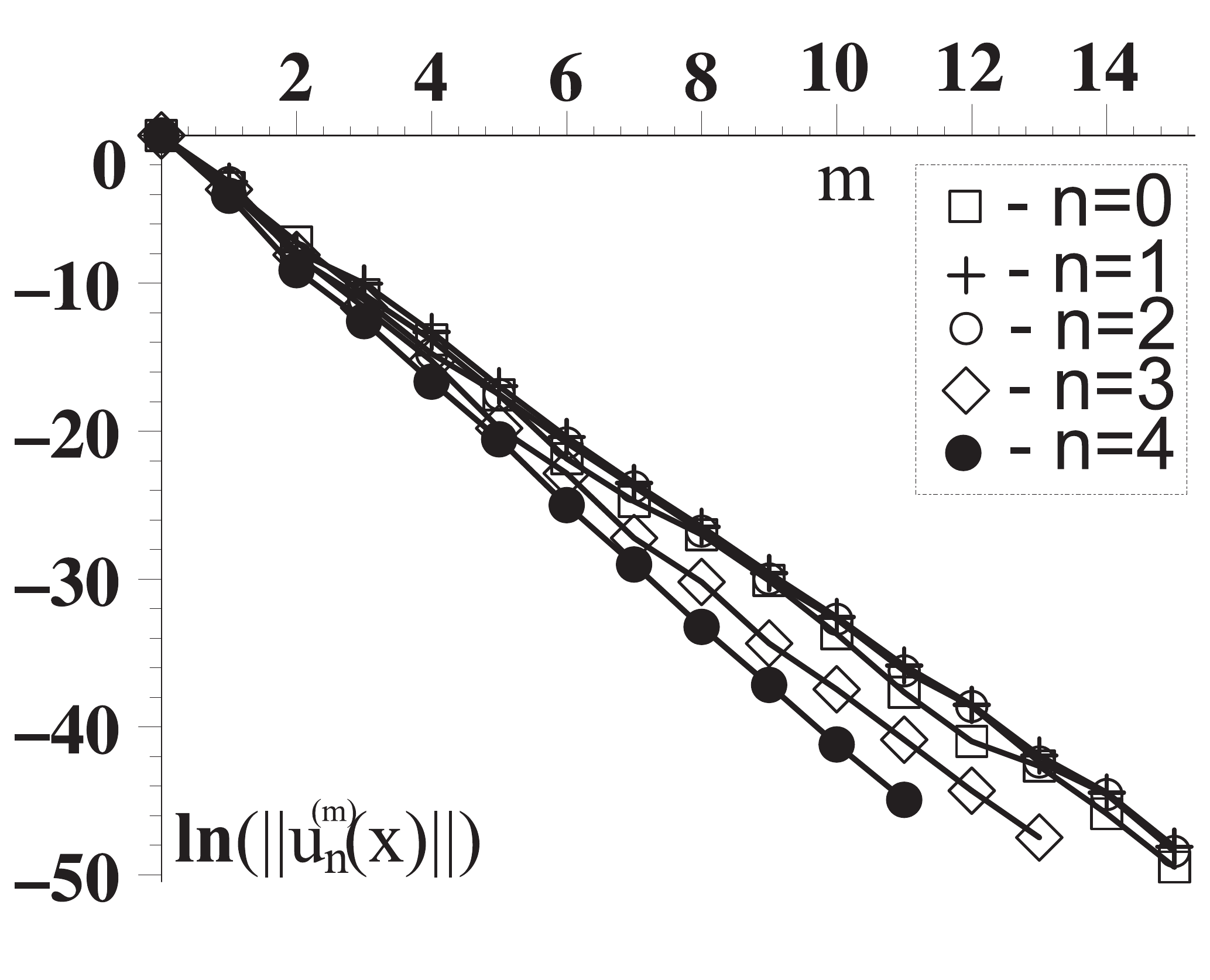}} \\ }
%\center{\includegraphics[width=1\linewidth]{r} \\ ?)}
%\caption{$\left|\nu_{2}\left(x\right)-u^{\ast}\left(x\right)\right|$}
\end{minipage}
\hfill
\begin{minipage}[h]{0.48\linewidth}
\center{\rotatebox{-0}{\includegraphics[height=0.8\linewidth,
width=1.0\linewidth]{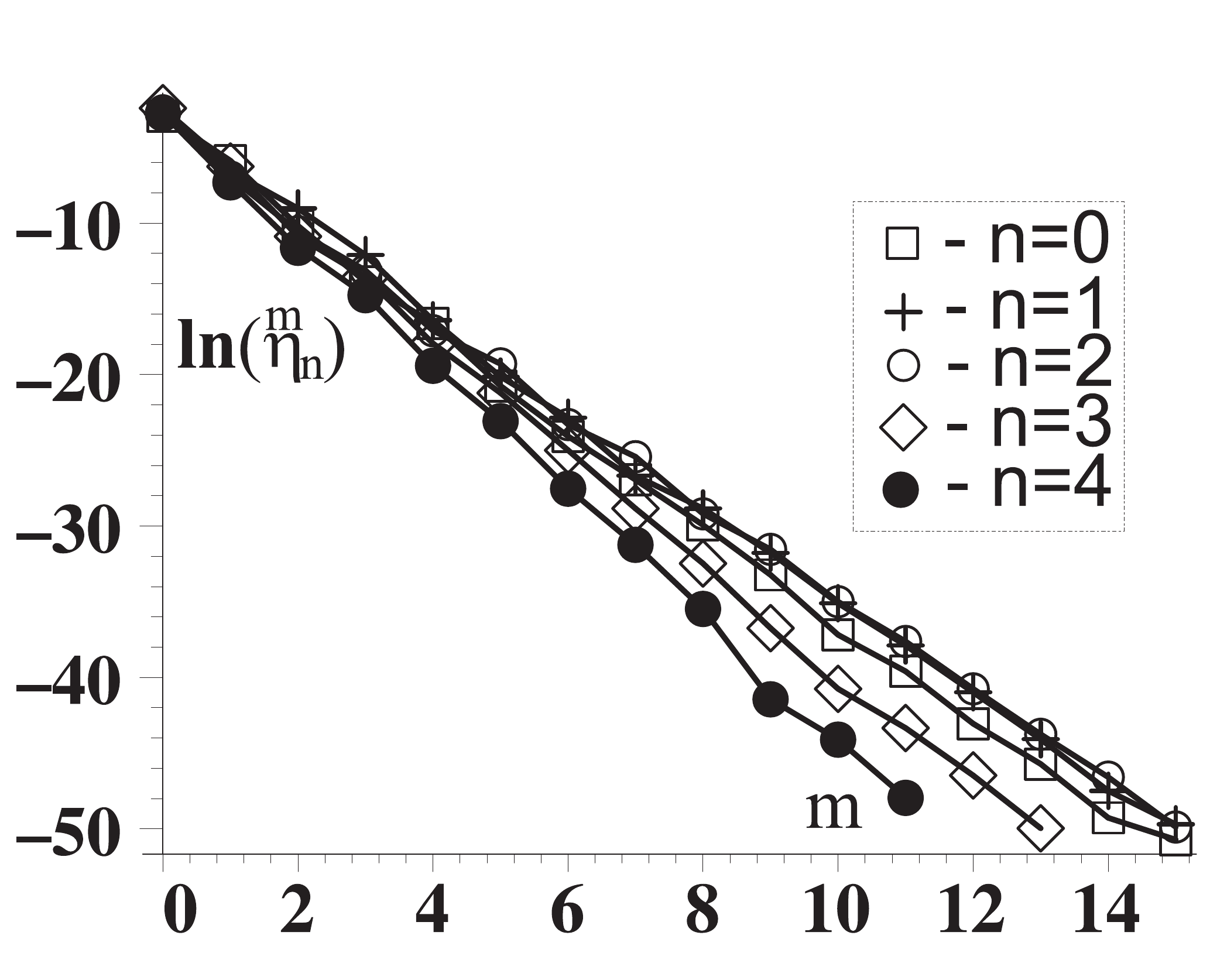}} }
%\center{\includegraphics[width=1\linewidth]{e} \\ ?)}
%\caption{$\left|\nu_{0}\left(x\right)-u^{\ast}\left(x\right)\right|$}
\end{minipage}
\caption{Example 1. The graphs of the functions $\ln\left(\left\|u^{(m)}_{n}(x)\right\|\right)$ (left)  and $\ln\left(\stackrel{m}{\eta}_{n}\right)$  (right).}\label{image1}
\end{minipage}
\end{figure}

Analyzing the data in tables \ref{table_2}, \ref{table_3} and \ref{table_4}, we can conclude that the FD-method with $N=1$ converges much more slower then the FD-method with $N=3.$ This fact is in good agreement with the results of Theorem \ref{main_theorem}. Furthermore, as it follows from table \ref{table_4}, the convergence rate of the FD-method increases as the index $n$ of the trial eigenvalue increases. Figure \ref{image1} illustrates the exponential nature of the FD-method's convergence.

{\bf Example 2.}
As the second example, let us consider problem \eqref{examp_1} with
\begin{equation}\label{examp_6}
    q(x)=\ln\left(\left|\left(\frac{5}{12}-x\right)\left(\frac{1}{3}+x\right)\right|\right).
\end{equation}
The results obtained with SLEIGN2 are presented in table \ref{table_5} below.
\begin{table}[htp]
\caption{Example 2, the data obtained by SLEIGN2}
\centering
\begin{tabular}{|c|c|c|c|}
  \hline
  $n$ & $\lambda_{n, sl2}$ & $TOL$ & $IFLAG$ \\
  \hline
  0 & $-1.98326983D+00$ & $0.46748D-08$ & $1$ \\
  \hline
  1 & $0.855187683D+00$ & $0.73426D-07$ & $1$ \\
  \hline
  2 & $0.489606686D+01$ & $0.35447D-07$ & $1$ \\
  \hline
  3 & $0.104183770D+02$ & $0.40228D-07$ & $1$ \\
  \hline
  4 & $0.188163965D+02$ & $0.61329D-11$ & $1$ \\
  \hline
\end{tabular}
\label{table_5}
\end{table}

\begin{figure}[htbp]
\begin{minipage}[h]{1\linewidth}
\begin{minipage}[h]{0.48\linewidth}
\center{\rotatebox{-0}{\includegraphics[height=0.8\linewidth,
width=1.0\linewidth]{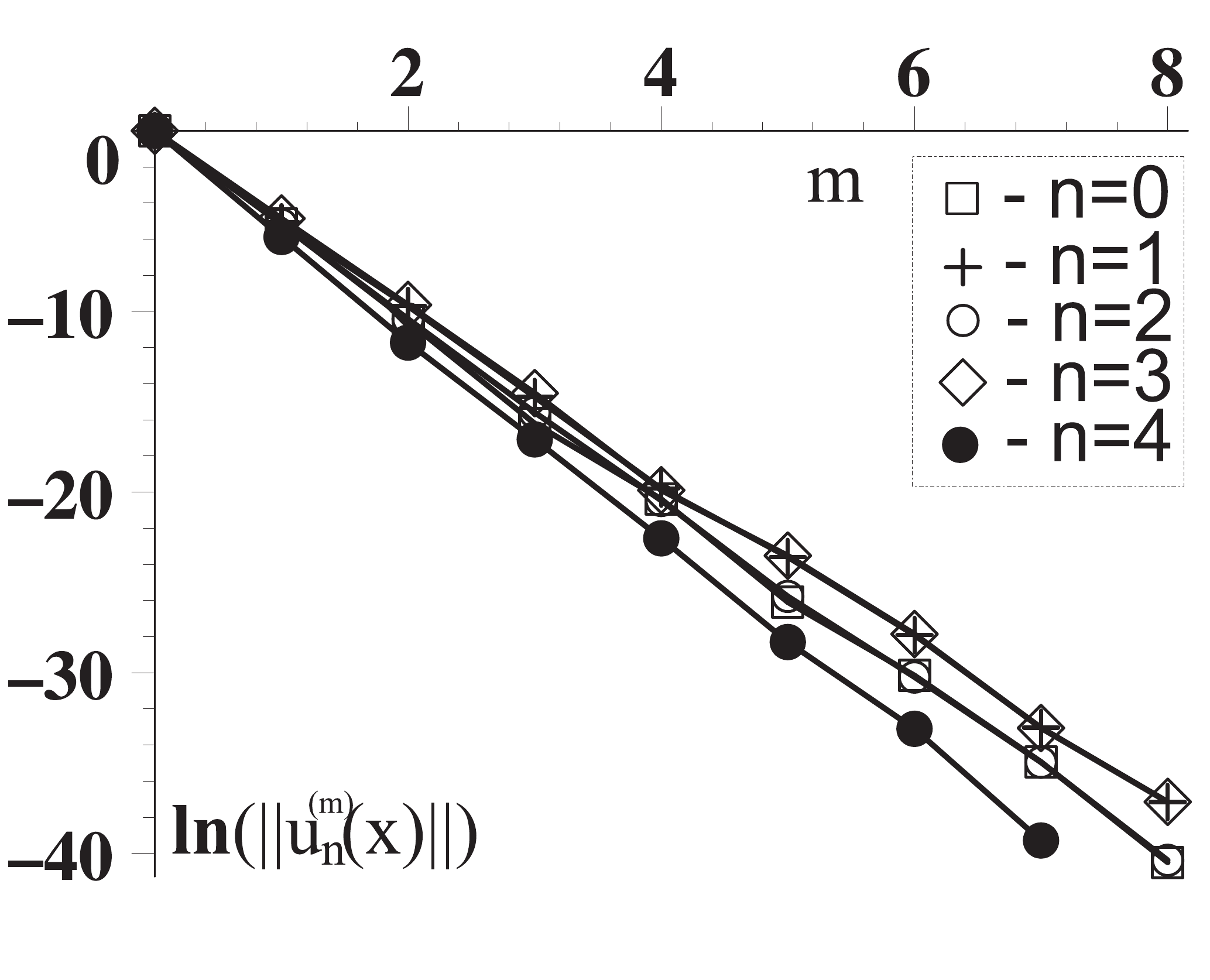}} \\ }
%\center{\includegraphics[width=1\linewidth]{r} \\ ?)}
%\caption{$\left|\nu_{2}\left(x\right)-u^{\ast}\left(x\right)\right|$}
\end{minipage}
\hfill
\begin{minipage}[h]{0.48\linewidth}
\center{\rotatebox{-0}{\includegraphics[height=0.8\linewidth,
width=1.0\linewidth]{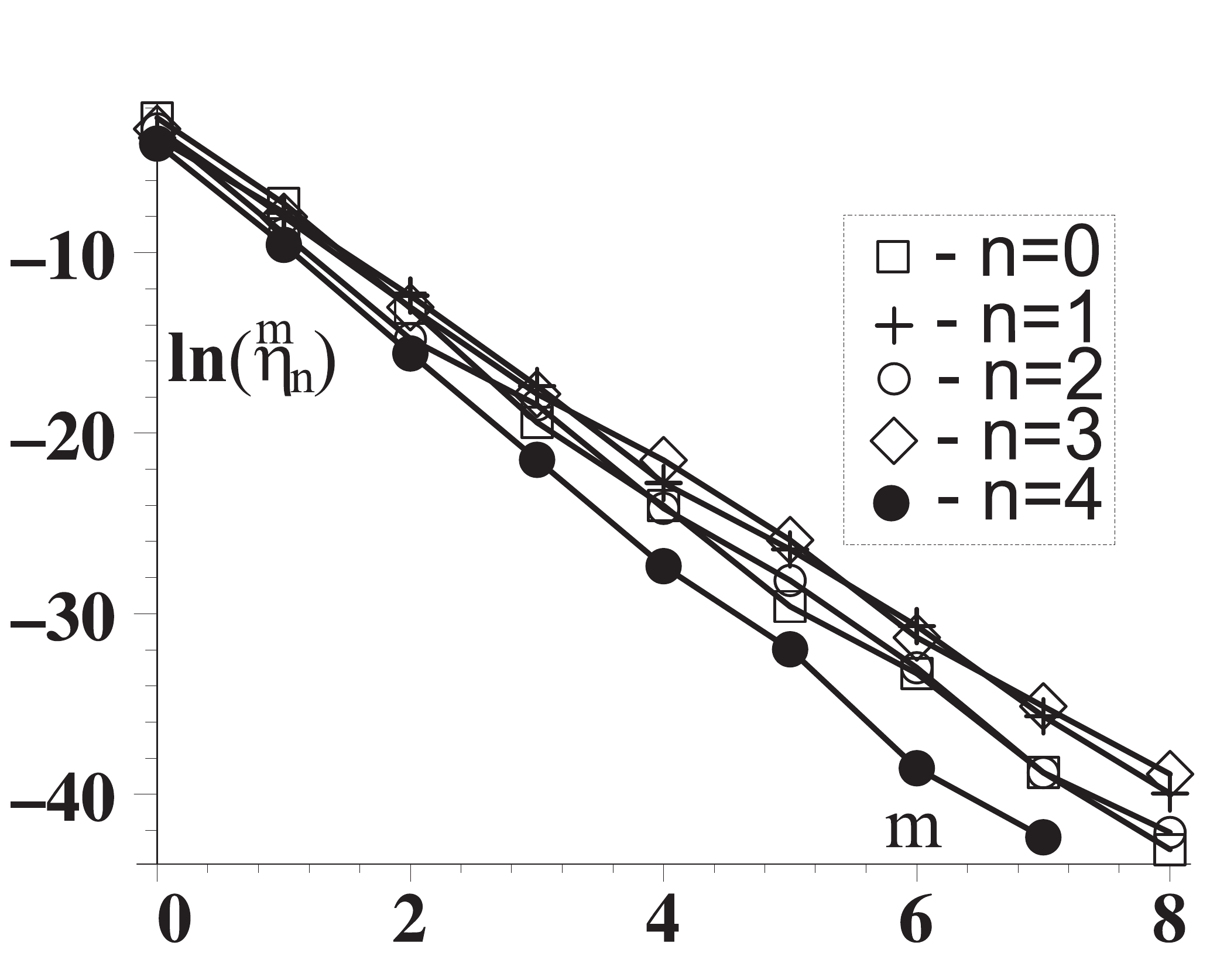}} }
%\center{\includegraphics[width=1\linewidth]{e} \\ ?)}
%\caption{$\left|\nu_{0}\left(x\right)-u^{\ast}\left(x\right)\right|$}
\end{minipage}
\caption{Example 2. The graphs of the functions $\ln\left(\left\|u^{(m)}_{n}(x)\right\|\right)$ (left)  and $\ln\left(\stackrel {m}{\eta}_{n}\right)$ (right).}\label{image2}
\end{minipage}
\end{figure}

It is worth to emphasize that the problem under consideration do not satisfy the conditions of theorem \ref{main_theorem}. However, from the results presented in table \ref{table_6} and figure \ref{image2} it follows that the FD-method have successfully handled this problem as opposed to SLEIGN2, which gives results with essential errors (see the leftmost colon in the table \ref{table_6}).

\begin{table}[htbp]
\caption{Example 2, FD-method with $N=24.$}
\centering
\begin{tabular}{|c|c|c|c|c|c|c|}
  \hline
  $n$ & $m$ & $\stackrel {m}{\lambda}_{n}$ & $\left|\lambda_{n}^{(m)}\right|$ & $\left\|u_{n}^{(m)}(x)\right\|$ & $\stackrel {m}{\eta}_{n}$ & $\left|\stackrel{m}{\lambda}_{n}-\lambda_{n, sl2}\right|$  \\
  \hline
  0 & 8 & $-1.9831442709774408386$ &$0.99\times 10^{-17}$& $0.26\times 10^{-17}$ & $0.20\times 10^{-18}$ & $0.125559\times 10^{-3}$\\
  \hline
  1 & 8 & $0.85727032837311800023$ &$0.47\times 10^{-15}$& $0.73\times 10^{-16}$ & $0.44\times 10^{-17}$ & $0.20826454\times 10^{-2}$\\
  \hline
  2 & 8 & $4.8939506826799075597$ &$0.98\times 10^{-17}$& $0.29\times 10^{-17}$ & $0.51\times 10^{-18}$ & $0.2116177\times 10^{-2}$\\
  \hline
  3 & 8 & $10.420511296257433545$ &$0.3\times 10^{-17}$& $0.74\times 10^{-16}$ & $0.13\times 10^{-16}$ & $0.213430\times 10^{-2}$\\
  \hline
  4 & 7 & $18.816396521508987920$ &$0.11\times 10^{-16}$& $0.87\times 10^{-17}$ & $0.39\times 10^{-18}$ & $0.2\times 10^{-7}$ \\
  \hline
\end{tabular}
\label{table_6}
\end{table}

\begin{figure}[htbp]
\begin{minipage}[h]{1\linewidth}
\begin{minipage}[h]{0.48\linewidth}
\center{\rotatebox{-0}{\includegraphics[height=0.8\linewidth,
width=1.0\linewidth]{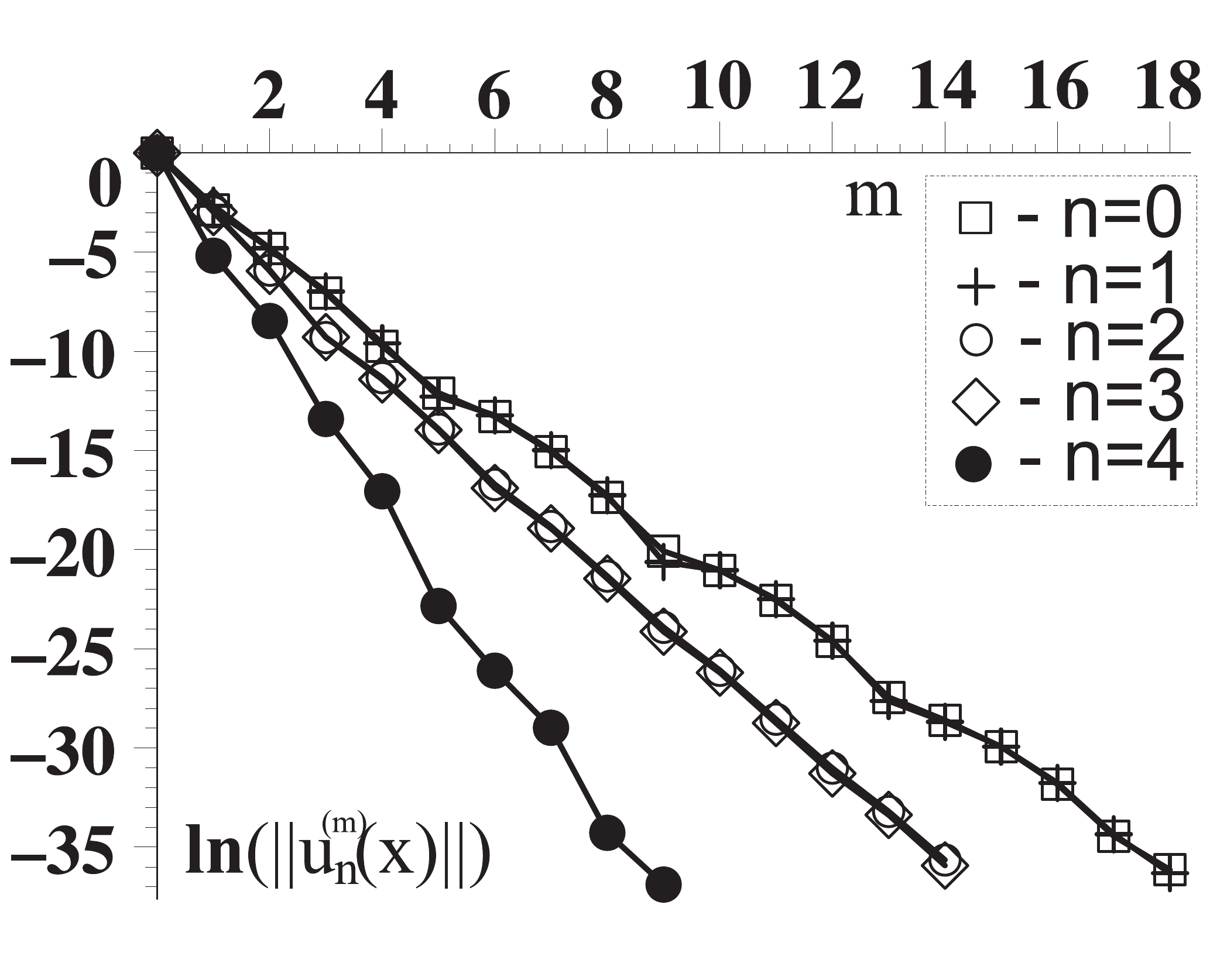}} \\ }
%\center{\includegraphics[width=1\linewidth]{r} \\ ?)}
%\caption{$\left|\nu_{2}\left(x\right)-u^{\ast}\left(x\right)\right|$}
\end{minipage}
\hfill
\begin{minipage}[h]{0.48\linewidth}
\center{\rotatebox{-0}{\includegraphics[height=0.8\linewidth,
width=1.0\linewidth]{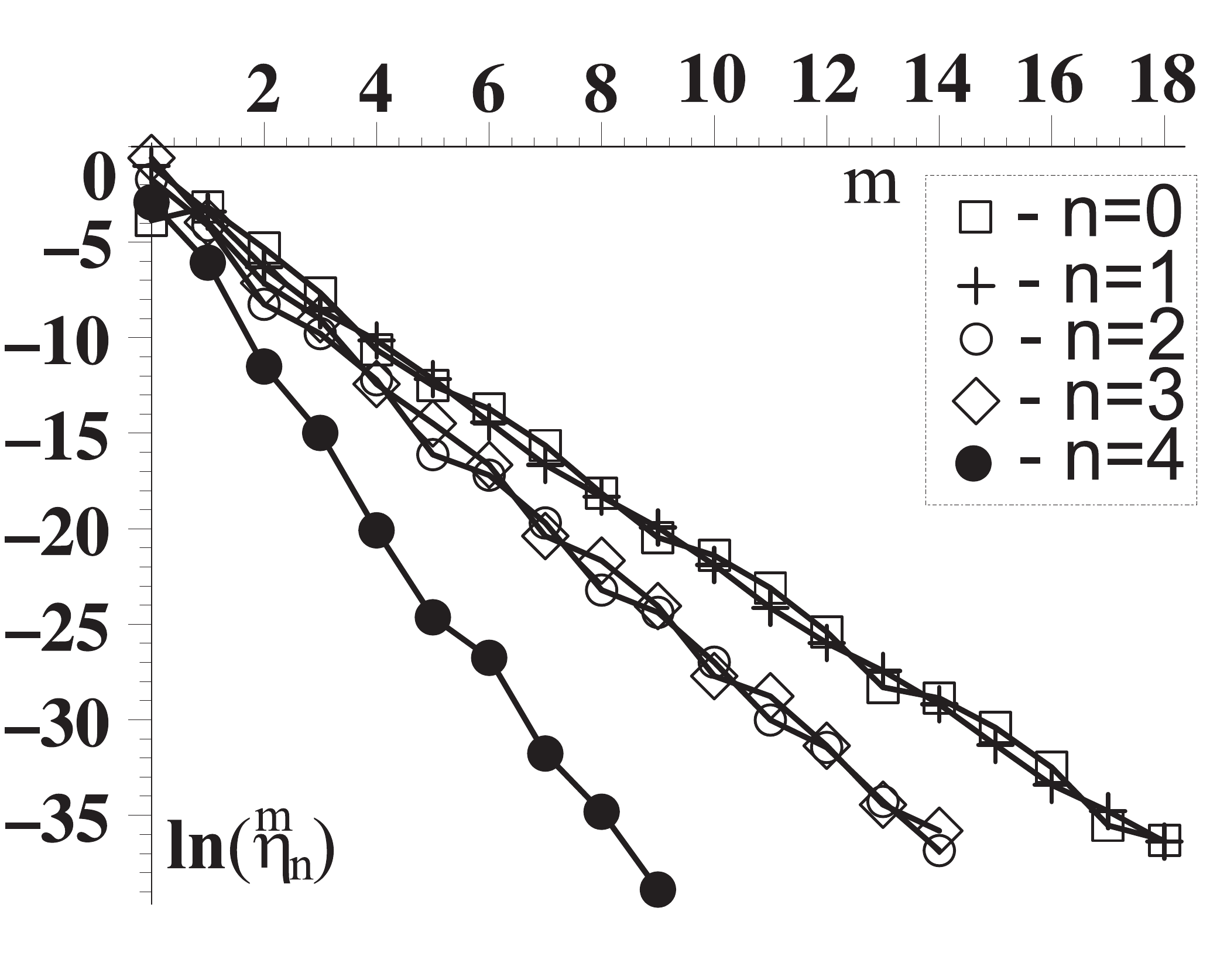}} }
%\center{\includegraphics[width=1\linewidth]{e} \\ ?)}
%\caption{$\left|\nu_{0}\left(x\right)-u^{\ast}\left(x\right)\right|$}
\end{minipage}
\caption{Example 3. The graphs of the functions $\ln\left(\left\|u^{(m)}_{n}(x)\right\|\right)$ (left)  and $\ln\left(\stackrel{m}{\eta}_{n}\right)$ (right).}\label{image3}
\end{minipage}
\end{figure}

{\bf Example 3.}
Finally, let us consider problem \eqref{examp_1} with
\begin{equation}\label{examp_3}
    q\left(x\right)=\frac{1}{\sqrt{\left|x+\frac{1}{3}\right|}}+\ln\left(\left|x-\frac{1}{3}\right|\right).
\end{equation}
We find out that SLAIGN2 does not handle this problem. But the FD-method does. The results obtained by the FD-method are presented in table \ref{table_7} and figure \ref{image3}. As before, figure \ref{image3} confirms the exponential convergence rate of the method.

\begin{table}[htbp]
\caption{Example 3, FD-method with $N=12.$}
\centering
\begin{tabular}{|c|c|c|c|c|c|}
  \hline
   $n$& $m$ & $\stackrel{m}{\lambda}_{n}$ &$\left|\lambda_{n}^{(m)}\right|$& $\left\|u_{n}^{(m)}(x)\right\|$ & $\stackrel{m}{\eta}_{n}$  \\
  \hline
   0 & 18 & $0.40796999146419634$ & $0.42\times 10^{-15}$ & $0.20\times 10^{-15}$ & $0.17\times 10^{-15}$ \\
  \hline
  1 & 18 & $3.4136861164474333$ &$0.4\times 10^{-15}$& $0.17\times 10^{-15}$ & $0.16\times 10^{-15}$ \\
  \hline
  2 & 14 & $6.7759537951814352$ & $0.2\times 10^{-15}$ & $0.33\times 10^{-15}$ & $0.10\times 10^{-15}$ \\
  \hline
  3 & 14 & $13.323487340142488$ & $0.2\times 10^{-14}$ & $0.25\times 10^{-15}$ & $0.28\times 10^{-15}$ \\
  \hline
  4 & 9 & $20.8431972121837340$ &$0.1\times 10^{15}$& $0.96\times 10^{-16}$ & $0.13\times 10^{-16}$ \\
  \hline
\end{tabular}
\label{table_7}
\end{table}

\section{Conclusions}
In the present paper we construct and theoretically justified the generalized algorithm of the FD-method for solving the Sturm-Liouville problem for differential equation of the second order \eqref{intro_exp_1}, \eqref{intro_exp_2} with piecewise continuous  functional coefficient $q(x).$ As it follows from Theorem \ref{main_theorem}, the generalized FD-method, which uses the piecewise constant approximation of the function $q(x),$ can be applied for the approximation of eigenvalues and eigenfunctions with any nonnegative index $n.$ The convergence rate of the method can be increased by decreasing the value $\left\|q(x)-\bar{q}(x)\right\|_{\infty, [-1,1]}.$ In the case when $\bar{q}(x)\equiv 0$ (this case was considered in \cite{Klimenko_Ya}) the FD-method (if converges) allows us to calculate the approximation to the eigensolution  analytically. But in general case, when $\bar{q}(x)\neq 0,$ the analytical calculations are almost always impossible and it is necessary to use numerical integration methods, such as sinc quadratures and Stenger's formula (see \cite{Kahaner_Nesh}, \cite{Stenger_num_meth}).

The problems considered in examples 2 and 3 do not satisfy the conditions of Theorem \ref{teor_1} because of the unboundedness of the function $q(x)$ on $[-1,1].$ However, as it was shown in the mentioned examples, the method successfully converges whereas the well known in the mathematical world package SLAIGN2 either gives not more then three correct numbers after decimal point (example 2) or can not handle the problem at all (example 3). This examples indicate that the FD-method has a considerable potential which are to be investigated in further mathematical works.

\bibliography{references_stat}
\end{document}